\documentclass[12pt]{article}
\usepackage{theorem,amsfonts}

\newtheorem{theorem}{Theorem}
\newtheorem{proposition}{Proposition}
\newtheorem{lemma}{Lemma}
\newtheorem{corollary}{Corollary}
\theorembodyfont{\upshape}

\newtheorem{remark}{Remark}
\newtheorem{proof}{Proof}

\newtheorem{acknowledgement}{Acknowledgement}

\newcommand{\bt}{\begin{theorem}}
\newcommand{\et}{\end{theorem}}
\newcommand{\bl}{\begin{lemma}}
\newcommand{\el}{\end{lemma}}
\newcommand{\bp}{\begin{proposition}}
\newcommand{\ep}{\end{proposition}}
\newcommand{\bo}{\begin{proof}}
\newcommand{\eo}{\end{proof}}
\newcommand{\br}{\begin{remark}}
\newcommand{\er}{\end{remark}}
\newcommand{\bc}{\begin{corollary}}
\newcommand{\ec}{\end{corollary}}
\newcommand{\be}{\begin{enumerate}}
\newcommand{\ee}{\end{enumerate}}

\title{Ergodic amenable actions of algebraic groups}
\author{C. R. E. Raja}
\date{}

\begin{document}
\maketitle

\let\epsi=\epsilon
\let\vepsi=\varepsilon
\let\lam=\lambda
\let\Lam=\Lambda 
\let\ap=\alpha
\let\vp=\varphi
\let\ra=\rightarrow
\let\Ra=\Rightarrow 
\let\LRa=\Leftrightarrow
\let\Llra=\Longleftrightarrow
\let\Lla=\Longleftarrow
\let\lra=\longrightarrow
\let\Lra=\Longrightarrow
\let\ba=\beta
\let\ga=\gamma
\let\Ga=\Gamma
\let\un=\upsilon

\begin{abstract}
We prove that every ergodic amenable action of an algebraic group over a 
local field of characteristic zero is induced from an ergodic action of an 
amenable subgroup.  
\end{abstract}


It was shown by Zimmer that an ergodic amenable action of a connected 
locally compact group is induced by an ergodic action of an 
amenable subgroup, here we prove the following analogue for algebraic 
groups over any local field of characteristic zero.  Let $\mathbb K$ be 
denote a local field of characteristic zero.  

\bt\label{mt}
Let $G$ be the group of $\mathbb K$-points of a Zariski-connected
algebraic group defined over $\mathbb K$.  Let $S$ be an ergodic amenable 
$G$-space.  Then the action on $S$ is induced from an ergodic action of 
an amenable subgroup of $G$.   
\et

In order to prove this theorem we first prove the analogue of
Moore's result (see 3.2.22 and 9.2.5 of \cite{Z}) and the rest is
similar to the proof of Zimmer.  Our approach is different from \cite{M} 
and we use a lemma of Furstenberg.  Let $V$ be a finite-dimensional vector 
space over $\mathbb K$.  Let $P(V)$ be the corresponding 
projective space.  Let $\Pi \colon V\setminus (0) \ra P(V)$ 
be the natural quotient map.  It is a well known fact that any element 
$g$ of $GL(V)$ gives a homeomorphism $\Pi (g)$ of $P (V)$.  Let $PGL(V)$
denotes the group consisting of $\Pi (g)$ for all $g\in GL(V)$: elements
of $PGL(V)$ are known as projective linear transformations.  
A subset $L$ of $V$ is called a {\bf \it quasi-linear variety} if it is a
union of finitely many subspaces of $V$.  We now state the following lemma 
due to Furstenberg (see \cite{F}).

\bl\label{fl}
Let $(g_n)$ be a sequence in $PGL (V)$ and $\mu$ and $\lam$ be probability 
measures on $P(V)$.  Then there is a subsequence $(g_{k_n})$ and a 
transformation $\tau$ of $P(V)$ onto $\Pi (L\setminus (0))$ where $L$ is a 
quasi-linear variety with the following properties: 

\be
\item $(g_{k_n})$ converges to $\tau$ pointwise on $P(V)$ and 
if $(g_{k_n})$ does not converge in $PGL(V)$, then $L$ is
a proper quasi-linear variety of $V$;

\item If $g_n (\mu )\ra \lam$ in the space of probability measures 
on $P(V)$, equipped with weak* topology, then $\lam$ is 
supported on $\Pi (L\setminus (0))$.
\ee
\el 

Let $X$ be any locally compact topological space and ${\cal P}(X)$ be the 
space of all regular Borel probability measures on $X$ equipped with the 
weak* topology which is the weakest topology on ${\cal P}(X)$ for which 
the functions $\mu \mapsto \mu (f)$ are continuous for all continuous 
bounded functions on $X$.  Let $G$ be any locally
compact group acting on $X$ by homeomorphisms of $X$.  Then the action of 
$G$ induces an action of $G$ on ${\cal P}(X)$.  
For a probability measure $\mu$ on $X$, we define the subgroups 
${\cal I}_G (\mu )\, =\, \{g\in G \, \mid \, g\mu \, =\, \mu \}$ and 
$I_G(\mu ) \,=\, \{\ap \in A\, \mid \, \ap (x) \, =\, x \, \, {\rm 
for\, all }\, \, x \in S(\mu )\}$ where $S(\mu )$ denotes the 
support of $\mu$.

\bp\label{dp}
Let $G$ be an algebraic  subgroup of $PGL (V)$ and $\mu$ be a
probability measure on $P(V)$.  Then ${\cal I}_G(\mu )/I_G(\mu )$ is a
compact group.  
\ep

\bo
Let $L$ be the smallest quasi-linear variety such that 
$\Pi (L\setminus (0))$ contains the support of $\mu$.  Let $H$ be the 
algebraic subgroup of $G$ consisting of transformations that preserve 
$\Pi (L\setminus (0))$, that is, 
$H=\{ g\in G\mid g(\Pi (L\setminus (0)))=\Pi (L\setminus (0))\}$.  
Then it is easy to see that ${\cal I}_G(\mu)$ 
is contained in $H$.  Let $(g_n)$ be a sequence in ${\cal I}_G(\mu)$.  
Then by Lemma \ref{fl} and by passing to a subsequence, if necessary, we 
may assume that there is transformation $\tau $ of $P(V)$ onto 
$\Pi (U\setminus (0))$ where $U$ is a quasi-linear variety such that 
$\Pi (U\setminus (0))$ contains the support of $\mu$ and $g_n \ra \tau$
pointwise on $P(V)$.  Since $g_n (\mu) \, =\, \mu$ for all $n \geq 1$,
$\tau (\mu)= \mu$ and hence for $n \geq 1$, 
$g_n (\Pi (L\setminus (0))) = \tau (\Pi (L\setminus (0)))\, 
=\Pi (L\setminus (0))$.  Let $L = \cup _{i=1}^k W_i$ where each $W_i$ is 
a subspace of $V$.  Then $H$ has a subgroup 
$N= \{g \in G \mid  g(\Pi (W_i \setminus (0)))= 
\Pi (W_i \setminus (0)) {\rm ~~for ~~ all}~~1\leq i \leq k\}$ of finite
index.  Now by passing to subsequence, if necessary, we may assume there
is a $h\in G$ such that $g_n h\in N$ for all $n \geq 1$.  Let $\ap$ be the 
restriction of $\tau h$ to $\Pi (L\setminus (0))$.  Then for 
$1\leq i\leq k$, 
$\ap (\Pi (W_i\setminus (0)))=\Pi (W_i\setminus (0))$.  
Let $\ap _i$ be the restriction of $\ap$ to $\Pi (W_i\setminus (0))$ for 
$1\leq i \leq k$.  
Let $\Phi \colon N \ra \prod _{i=1}^k PGL(W_i)$ be defined by 
$\Phi (g) =(\Phi _i (g))_{i=1}^k$ where $\Phi _i(g)$ is the restriction of 
$g$ to $\Pi (W_i\setminus (0))$ for $1\leq i \leq k$.  
Since $N$ is an algebraic group, the image $\Phi (N)$ is closed and it is 
isomorphic to $N/{\rm ker\Phi}$.  
Since $\Phi (g_nh) \ra (\ap _i)_{i=1}^k$ and the 
kernel of $\Phi$ is $I_G(\mu )$, $(g_nh)$ is relatively compact in 
$N/I_G(\mu )$.  This proves the proposition since ${\cal I}_G(\mu )$ is a 
closed subgroup of $H$ containing $I_G(\mu )$. 
\eo

By an algebraic group over $\mathbb K$ we mean the group of 
$\mathbb K$-points of an algebraic group defined over $\mathbb K$.

\bc\label{dc}
Let $G$ be an algebraic group over $\mathbb K$ and 
$H$ be an algebraic subgroup of $G$.  Let $\mu$ be a 
probability measure on $G/H$.  Then $G$ acts on $G/H$ in a canonical way.
Let ${\cal I}(\mu)= \{g\in G \mid  g \, \, {\rm action \, \, on}
\,\, G/H \, \, {\rm preserves}\, \, \mu \}$ and $I(\mu ) = \{g\in G
\mid  g \, \, {\rm acts \, \, trivially\, \, on \, \, support \, \,
of}\, \, \mu \}$.  Then ${\cal I}(\mu )/I(\mu )$ is a compact group.
\ec

\bo
Since $H$ is an algebraic subgroup of an algebraic group $G$, there is a 
finite-dimensional vector space $V$ on which $G$ acts by linear 
transformations and there is a vector $v_0\in V$ such that 
$H = \{ g\in G \mid  gv_0\in (v_0) \}$ where $(v_0)$ is the 
one-dimensional subspace of $V$ spanned by $v_0$.  
This implies that $G/H$ can be viewed as a Borel subset (in fact, a
locally closed subset) of $P(V)$ 
and $G$ acts on $P (V)$ as an algebraic group.  Now the result 
follows from Proposition \ref{dp}.      
\eo  

\bc\label{dp1}
Let $G$ be the group of $\mathbb K$-points of a Zariski-connected
algebraic group defined over $\mathbb K$.  Let $P$ be a minimal parabolic
subgroup of $G$.  Then for any algebraic subgroup $H$ of $G$ and any 
$\mu \in {\cal P}(G/P)$, ${\cal I}_H(\mu )=\{h\in H \mid h\mu = \mu \}$
is amenable.  In particular, a subgroup of $G$ is amenable if and only if
it has a fixed point in the space of probability measures on $G/P$.
\ec

\bo
By Corollary \ref{dc}, ${\cal I}_H(\mu )/I_H(\mu )$ is
compact where $I_H(\mu )= \{h\in H \mid  hx= x\, \, 
{\rm for \, \, all \, \, }x {\rm \, \, in \, \, support \, \, of \, \, }
\mu \}$.  Hence it is enough to prove $I_H(\mu )$ is amenable.  Since
elements of $I_H(\mu )$ fixes the support of $\mu$ pointwise, $I_H(\mu )$
is contained in a conjugate of $P$ and hence since $P$ is amenable 
$I_H(\mu )$ is amenable.     
\eo

\bo ({\bf Proof of Theorem \ref{mt}})
Consider the cocycle $\ap$ defined by $\ap (s, g) = g$ for all $s \in S$
and $g \in G$.  Let $N$ be the solvable radical of $G$.  Then $H= G/N$ is
semisimple and let $P$ be the minimal parabolic subgroup of $H$.  By
amenability there is a $\ap$-invariant function $f \colon S \ra {\cal
P}(H/P)$.  Since the action of $G$ on ${\cal P}(H/P)$ is smooth (see
Corollary 3.2.17 of \cite{Z}) and by Corollary \ref{dp1}, stabilizers are
amenable.  By cocycle reduction Lemma 5.2.11 of \cite{Z}, 
$\ap$ is equivalent to a cocycle taking values in an amenable subgroup, 
say, $M$ of $G$.  We now claim that the action of $G$ on $S$ is induced 
from action of $M$.  By 4.2.18 of \cite{Z}, there is an 
$\ap$-invariant function $\phi \colon S\ra G/M$.  This implies by Theorem 
2.5 or Corollary 2.6 of \cite{Z0} that the action on $S$ is induced from 
an ergodic amenable (because $M$ is amenable) action of $M$.  
\eo

It may be recalled that any local field $\mathbb K$ of characteristic 
zero is either $\mathbb R$ or $\mathbb C$ or finite extension of a 
$p$-adic field.  Zimmer's result Theorem 5.7 of \cite{Z0} covers the case 
of real algebraic groups and all connected groups (see 9.2.5 of \cite{Z}).  
We denote by ${\mathbb Q}_p$ the $p$-adic field.  Let $\tilde G$ be an 
algebraic group defined over ${\mathbb Q}_p$ and $G= \tilde 
G({\mathbb Q}_p)$, the group of ${\mathbb Q}_p$-points of $\tilde G$.  
Then it may be easily seen that $G$ is a totally disconnected group, hence 
these groups are not covered by Zimmer's result.  Now we will present a 
few explicit examples of these groups.  
\be

\item The elementary examples are the additive group ${\mathbb Q}_p$ and 
the multiplicative group ${\mathbb Q}_p\setminus (0)$ known as the 
one-dimensional split torus.

\item $GL_n ( {\mathbb Q}_p)$ be the group of all invertible $n\times n$ 
matrices with entries in ${\mathbb Q}_p$, in general, all invertible 
elements in a finite-dimensional algebra over ${\mathbb Q}_p$.

\item $SL_n( {\mathbb Q}_p)= \{A \in GL_n( {\mathbb Q}_p) \mid \det (A) 
=1 \}$, the special linear group.

\item For any $m \in {\mathbb N}$ and any $2m\times 2m$ skew-symmetric 
matrix $E$, that is $E^t = -E$, define the symplectic group 
$SP_{2m}(E) = \{ A \in GL_{2m}({\mathbb Q}_p) \mid A^tEA=E \}$.  

\item The group of upper triangular matrices, $UT_n({\mathbb Q}_p) = \{ 
(g_{ij}) \in GL_n({\mathbb Q}_p) \mid g_{ij} = 0 ~~ {\rm if }~~ j<i \}$ 
and the group of all unipotent martices, $U_n({\mathbb Q}_p) = 
\{ (g_{ij}) \in UT_n({\mathbb Q}_p) \mid g_{ii} = 1 \}$.  

\ee
These groups and their direct and semi-direct products (in some cases) 
and many other classical algebraic groups are covered by our result but 
not by Zimmer's result.  Though the groups in examples $1$ and $5$ are 
amenable but their direct and semidirect products with $GL_n$ and $SL_n$ 
or some suitable subgroups of $GL_n$ and $SL_n$ are not amenable, 
for example, the general affine group 
${\mathbb Q}_p^n \times _s GL_n $ is not amenable.  

We will end the article with few consequences of the main result.  In
\cite{Z0}, it is shown that for any lattice $\Gamma$ in 
$SL(2, {\mathbb C})$, the $\Gamma$-space $SL(2, {\mathbb C})/ N$ 
(where $N$ is the group of all 
upper triangular matrices in $SL(2,{\mathbb C})$) is ergodic amenable but 
it is not induced from an amenable subgroup.  In fact, Theorem \ref{mt} 
combined with arguments preceding Theorem 5.11 of \cite{Z0}, 
we conclude that any ergodic 
amenable $\Gamma$-space has a factor of the form $SL(2, {\mathbb C})/R$ 
for some amenable subgroup $R$ of $SL(2, {\mathbb C})$.  Also, results 
regarding minimal ergodic amenable actions of closed subgroups of 
$GL_n (\mathbb K)$ may also be obtained as in Theorem 5.9 and Corollary 
5.12 of \cite{Z0}.   

\begin{acknowledgement} 
I wish to thank the referee for useful remarks.  
\end{acknowledgement}

\vskip 0.25in

\noindent {C. Robinson Edward Raja, \\
Indian Statistical Institute, \\
Statistics and Mathematics Unit,\\
8th Mile Mysore Road,\\
R. V. College Post,\\
Bangalore - 560 059.\\
India.}

\noindent {creraja@isibang.ac.in}


\begin{thebibliography}{4}

\footnotesize


\bibitem {F} H. Furstenberg, A note on Borel density theorem, Proc. AMS
55 (1976), 209-212.

\bibitem {M} C. C. Moore, Amenable subgroups of semisimple groups and
proximal flows, Israel Journal of Mathematics 34 (1979), 121-138.   

\bibitem {Z0} R. J. Zimmer, Induced and amenable ergodic actions of Lie 
groups, Ann. Sci. École Norm. Sup. 11 (1978), no. 3, 407--428.

\bibitem {Z} R. J. Zimmer, {\it Ergodic Theory and Semisimple Groups}, 
Birkauser, Boston, (1984)

\end{thebibliography}
\end{document}